Nenad Antonić & Darko Mitrović & Tomislav Perić

# Orthogonality of H-distributions and applications

## Abstract

We extend Gérard's results on orthogonality of $L^2_{\text{loc}}$ sequences as a consequence of mutual singularity of corresponding H-measures (microlocal defect measures) to $L^p/L^q$ sequences and newly introduced notion of orhogonality for H-distributions. We apply the result to a homogenisation problem for the heterogeneous Boltzmann equation with space-dependent drift and periodic opacity.

Department of Mathematics
Faculty of Science
University of Zagreb
Bijenička cesta 30
Zagreb, Croatia

Faculty of Mathematics
University of Vienna
Oskar-Morgenstern-Platz 1
Vienna, Austria

on leave from:
Faculty of Maths and Science
University of Montenegro
Cetinjski put bb
Podgorica, Montenegro

Department of Mathematics
Faculty of Science
University of Zagreb
Bijenička cesta 30
Zagreb, Croatia

`nenad@math.hr`     `darko.mitrovic@univie.ac.at`     `tomislav.peric@math.hr`

This work is supported in part by the Croatian Science Foundation under projects , and by Croatian-Austrian bilateral project .

24<sup>th</sup> October 2025

# 1. Introduction

H-measures (or microlocal defect measures) were introduced some 35 years ago [14, 23] as mathematical objects capable of capturing limits of products of weakly converging sequences in $L^2$. Since then, the theory has been extended to H-distributions [4, 6, 7], in order to encompass products of $L^p$ and $L^q$ ($q \geqslant p'$, $1/p + 1/p' = 1$) sequences.

One of the original successes of H-measures was the connection of orthogonality of $L^2$ sequences and the measure-theoretic mutual singularity of the corresponding H-measures [14]. More precisely, we say that two weakly converging sequences $\mathsf{u}_n \rightharpoonup \mathsf{u}$ and $\mathsf{v}_n \rightharpoonup \mathsf{v}$ in $L^2_{\mathrm{loc}}(\Omega; \mathbf{C}^r)$ are *orthogonal* if we can pass to the limit in their product (by $\cdot$ we denote the standard sesquilinear scalar product on $\mathbf{C}^r$)

$$\mathsf{u}_n \cdot \mathsf{v}_n \overset{*}{\rightharpoonup} \mathsf{u} \cdot \mathsf{v} \qquad \text{in } \mathcal{D}'_0(\Omega) \,.$$

Of course, $\mathsf{u}_n \cdot \mathsf{v}_n$ is a bounded sequence in $L^1_{\mathrm{loc}}(\Omega; \mathbf{C})$, which has an accumulation point after embedding in the space of Radon measures (distributions of order zero) $\mathcal{D}'_0(\Omega) = (C_c(\Omega))'$, in the weak-$*$ topology. The orthogonality requires that the only accumulation point is $\mathsf{u} \cdot \mathsf{v}$. Note that the above weak-$*$ convergence is equivalent to weak-$*$ convergence in $\mathcal{D}'(\Omega) = (C_c^\infty(\Omega))'$, when restricted to $\mathcal{D}'_0(\Omega)$. However, $\mathcal{D}'(\Omega)$ being a reflexive space, it is equivalent to the weak convergence in $\mathcal{D}'(\Omega)$, which is the usual way of writing it in the literature.

A well-known limitation of H-measures is their $L^2$-character. In other words, the orthogonality of $(u_n)$ and $(v_n)$ can be characterised via H-measures only if both sequences are bounded in $L^p_{\mathrm{loc}}$ for $p \geqslant 2$ (so we essentially have the $L^2$ framework).

The aim of this paper is to go beyond the $L^2$-framework and to describe the orthogonality of $(u_n)$ and $(v_n)$ when the integrability of one of them is lower than $p = 2$. This will enable us to generalise the homogenisation result for heterogenenous Boltzmann type equations given in [14, Theorem 3.6] in the $L^p$, $p < 2$, setting, as well as the one from [13, Theorem 4.2] by taking local particle velocity (drift) in position–velocity space, which may depend on the environment (medium heterogeneity, refractive index, external forces); see (6) below. However, our result will be based solely on the generalisation of orthogonality to H-distributions (and their known properties), without any need to invoke the velocity averaging lemmas and drift non-degeneracy requirements.

We note that the equations considered in [13] and [14] contain additional terms, which are irrelevant in the framework of our approach and can easily be included.

As we are going to use the tool of H-distributions for that, we refer the reader to Appendix 1 where we state the main result on the existence of H-distributions [6], which also includes the result for H-measures.

Let us now briefly motivate the homogenisation problem on which we are going to apply our tools. Homogenisation of kinetic equations has been motivated by practical scenarios where transport coefficients oscillate rapidly across scales, such as the resonant neutron cross-sections in composite nuclear materials. Classical approaches replace the microscopically oscillatory linear Boltzmann equation with an effective model having smoother *homogenised* coefficients. The literature is rich in results under periodic or scale-separated assumptions—typically oscillations in one variable (often spatial)—and tools such as two-scale convergence or H-convergence are employed to handle the microstructure. For instance, simultaneous diffusion-homogenisation limits have been rigorously obtained for spatially heterogeneous linear Boltzmann models [8, 9, 15], while Dumas and Golse [13] showed in a periodic setting that the homogenised equation can retain the structure of the original transport model without introducing artificial memory. By contrast, homogenisation in the energy variable [17] produces an effective model with an explicit time-convolution kernel, exemplifying the appearance of memory effects. This phenomenon was first highlighted in Tartar's seminal work [21, 22] (cf. also [1]), which showed that even a simple ordinary differential equation with rapidly oscillating coefficients may converge to a nonlocal-in-time limit.

The central difficulty is caused by the fact that the homogenised equation may lose the semigroup (Markovian) property of the original dynamics. To overcome this, the *extended phase space method* was introduced [10, 11]: by augmenting the system with an auxiliary *age* variable, one derives a homogenised transport equation that preserves the local structure and generates a Markov semigroup in the extended variables. Only upon projection back to the physical phase space does one recover the integro-differential formulation with memory. This viewpoint reconciles the structural robustness of kinetic transport with the effective memory law observed in homogenised models.

As a final remark, in many applications the streaming operator is more general than the idealised term $\boldsymbol{\lambda} \cdot \nabla_{\mathbf{x}} f$. One often encounters $\partial_t f + \mathsf{a}(\mathbf{x}, \boldsymbol{\lambda}) \cdot \nabla_{\mathbf{x}} f$, where the drift field $\mathsf{a}(\mathbf{x}, \boldsymbol{\lambda})$ represents the local particle velocity. Physically, this occurs in radiative transfer through stellar atmospheres, where the photon group velocity depends on the refractive index [27], or in charged-particle transport under external electromagnetic fields. In neutron transport and related phenomena, strongly heterogeneous scattering and absorption cross-sections arise from the composite structure of the medium. In all such settings, direct numerical resolution of the microscopic Boltzmann equation is prohibitively complicated, and homogenisation provides effective macroscopic models that capture the essential physics while retaining computational tractability. The question of homogenisation for this type of problems has not yet been considered in the $\mathrm{L}^p$, $p < 2$ setting.

The paper is organised as follows. After the Introduction, in Section 2 we establish the notion of orthogonality for general distributions and recall the corresponding result for Radon measures. In Section 3 we extend this concept to the framework of H-distributions and formulate the main orthogonality theorem. Section 4 is devoted to the localisation principle, which generalises the compactness by compensation method to the setting of H-distributions. Finally, in Section 5 we present an application to the homogenisation of heterogeneous Boltzmann equations, where we illustrate how the abstract results developed earlier can be employed in a concrete kinetic setting.

## 2. Orthogonality of distributions

Let us now precisely state a simple orhogonality result [14, Proposition 3.1]:

*For two bounded sequences in $\mathrm{L}^2_{\mathrm{loc}}(\Omega)$ to be orthogonal it is enough that the corresponding H-measures $\mu$ and $\nu$ are (mutually) singular measures.*

Our goal is to extend this result to H-distributions (cf. Theorem 5). In order to do that, we first need to extend the notion of mutual singularity from Radon measures to distributions of finite order. Of course, we cannot use the ideas related to the Radon-Nikodým theorem as it is well known that such results for distributions of order more than zero are not valid.

In order to avoid ambiguity when the Riesz representation theorem is implicitly used, we shall reserve the term *mutually singular measures* for measures $\mu$ and $\nu$ considered as set functions, when it means that there is a measurable set $M$ such that $|\mu|(M) = |\nu|(cM) = 0$.

For real Radon measures, i.e. distributions of order zero, which form a Riesz space (i.e. lattice-ordered vector space) [12, II.1] we shall use the term *alien* (french: *étranger*), meaning that $\inf\{|\mu|, |\nu|\} = 0$. This definition extends to complex measures, as only the total variation of measures is used in the definition.

Additionally, for an arbitrary positive Radon measure $\mu$ the notion of upper integral $\mu^*$ can be defined for arbitrary positive functions (see [12, IV.3]). Let us emphasise that an arbitrary positive Radon measure and its upper integral coincide on the space of integrable positive functions.

Following [12, V.54], a complex Radon measure $\mu$ is *concentrated* on a subset $M$, if its complement $cM$ is locally negligible for $\mu$ (i.e. the intersection of $cM$ with any compact is negligible). The following proposition is valid [id., Proposition 13]:

*Two measures $\mu$ and $\nu$ are mutually alien if and only if there are two disjoint sets $M$ and $N$, such that $\mu$ is concentrated on $M$ and $\nu$ on $N$ (these sets may be taken to be universally measurable).*

This result equivalences alienness and certain orthogonality, which will not be feasible for extension to distributions of finite order. We shall slightly modify the above and baptise it as *orthogonality*, denoted by $\mu \perp \nu$. By $X$ we denote the domain, a differentiable manifold of sufficient smoothness.

Two distributions of finite order $\mu$ and $\nu$ with orders respectively $k$ and $l$ defined on $X$ are *orthogonal*, $\mu \perp \nu$, if

for any $\varepsilon > 0$ there are open sets $U$ and $V$ in $X$, $U \cup V = X$, such that

$$(\forall \varphi_1 \in \mathrm{C}_c^k(X))(\forall \varphi_2 \in \mathrm{C}_c^l(X)) \quad \operatorname{supp} \varphi_1 \subseteq U \ \& \ \operatorname{supp} \varphi_2 \subseteq V \Rightarrow \begin{array}{l} |\langle \mu, \varphi_1 \rangle| \leqslant \varepsilon \|\varphi_1\|_{\mathrm{C}_b^k(X)} \\ |\langle \nu, \varphi_2 \rangle| \leqslant \varepsilon \|\varphi_2\|_{\mathrm{C}_b^l(X)} \end{array}.$$

As a first step, let us prove that alienness and this orthogonality are equivalent for Radon measures ($k = l = 0$ in the definition).

**Theorem 1.** *Radon measures $\mu, \nu \in (\mathrm{C}_c(X))'$ are orthogonal if and only if they are mutually alien.*

Dem. Assuming alienness, the above equivalence provides two disjoint sets, $M \dot\cup N = X$, where $\mu$ is concentrated on $M$, and $\nu$ on $N$. As $\mu$ is negligible on $N$, for given $\varepsilon > 0$ we can find an open set $U \subseteq X$, such that $N \subseteq U$ and $|\mu|(U) < \varepsilon$. Similarly, we can find an open $V \subseteq X$ such that $M \subseteq V$ and $|\nu|(V) < \varepsilon$. Let $\varphi_1, \varphi_2 \in \mathrm{C}_c(X)$ such that $\operatorname{supp} \varphi_1 \subseteq U$ and $\operatorname{supp} \varphi_2 \subseteq V$. Now, it can easily be seen that

$$\begin{aligned} |\langle \mu, \varphi_1 \rangle| &\leqslant \langle |\mu|, |\varphi_1| \rangle \\ &\leqslant \langle |\mu|^*, \|\varphi_1\|_{\mathrm{C}_b(\Omega)} \chi_U \rangle \\ &= \|\varphi_1\|_{\mathrm{C}_b(\Omega)} |\mu|(U) \qquad \leqslant \varepsilon \|\varphi_1\|_{\mathrm{C}_b(X)}. \end{aligned}$$

Using the same argument we also get that

$$|\langle \nu, \varphi_2 \rangle| \leqslant \varepsilon \|\varphi_2\|_{\mathrm{C}_b(X)}.$$

Conversely, let us first note that if one of the measures is trivial, the measures are mutually alien, so the result trivially holds.

Let $\mu, \nu \in (\mathrm{C}_c(X))'$ be nontrivial Radon measures such that $\mu \perp \nu$. Furthermore, for $n \in \mathbf{N}$ by $U_n$ and $V_n$ we denote open sets $U$ and $V$ from the previous definition for $\varepsilon = 1/n$. Without loss of generality we can assume that $(U_n)$ and $(V_n)$ are decreasing sequences of sets, while $\bigcap_{n \in \mathbf{N}} U_n$ and $\bigcap_{n \in \mathbf{N}} V_n$ are nonempty sets. Define $R$ and $S$ to be their complements

$$R := X \setminus \bigcap_{n \in \mathbf{N}} U_n \qquad \text{and} \qquad S := X \setminus \bigcap_{n \in \mathbf{N}} V_n \ .$$

Since the characteristic function of an open set is lower semicontinuous [12, IV.4.2], it can easily be seen that

$$\begin{aligned} |\mu|^*(X \setminus R) &\leqslant |\mu|^*(U_n) \\ &= |\mu|^*(\chi_{U_n}) \\ &= \sup_{\substack{f \in \mathrm{C}_c(X) \\ 0 \leqslant f \leqslant \chi_{U_n}}} \langle |\mu|, f \rangle = \sup_{\substack{f \in \mathrm{C}_c(X) \\ 0 \leqslant f \leqslant \chi_{U_n}}} \sup_{\substack{g \in \mathrm{C}_c(X) \\ |g| \leqslant f}} |\langle \mu, g \rangle|. \end{aligned}$$

For arbitrary $f \in \mathrm{C}_c(X)$ such that $0 \leqslant f \leqslant \chi_{U_n}$ we have that $\operatorname{supp} f \subseteq U_n$. Therefore, we have the inequalities

$$\sup_{\substack{g \in \mathrm{C}_c(X) \\ |g| \leqslant f}} |\langle \mu, g \rangle| \leqslant \frac{1}{n} \sup_{\substack{g \in \mathrm{C}_c(X) \\ |g| \leqslant f}} \|g\|_{\mathrm{C}_b(X)} \leqslant \frac{1}{n}.$$

Now, it follows that $|\mu|^*(X \setminus R) \leqslant \frac{1}{n}$ for every $n \in \mathbf{N}$. Therefore, $|\mu|^*(X \setminus R) = 0$, which means that $\mu$ is concentrated on $R$. Similarly, we can get that $\nu$ is concentrated on $S$. Finally, we shall prove that $R$ and $S$ are disjoint. Let us suppose otherwise. Then there exist $x \in X$ and $n, m \in \mathbf{N}$ such that $x \notin U_n$ and $x \notin V_m$. Let $k = \max\{n, m\}$. Then, $x \notin U_k \cup V_k = X$ and we clearly arrive to a contradiction.

**Q.E.D.**

## 3. Result on orthogonality of H-distributions

From now on we shall consider $\Omega \subseteq \mathbf{R}^d$ to be a nonempty open set. We would like to apply general results of the previous section to H-distributions, which are known to be anisotropic distributions on the cospherical bundle $S^*\Omega = \Omega \times S^{d-1}$. Indeed, previous definitions can easily be generalised to anisotropic distributions (in particular on $S^*\Omega$, the distributions of order $m$ in $\mathbf{x} \in \Omega$ and order $n$ in $\boldsymbol{\xi} \in S^{d-1}$) by considering multiindices $\boldsymbol{\alpha} = (\boldsymbol{\alpha}_1, \boldsymbol{\alpha}_2)$ where $\boldsymbol{\alpha}_1 \in \mathbf{N}_0^d$ and $\boldsymbol{\alpha}_2 \in \mathbf{N}_0^{d-1}$ are such that $|\boldsymbol{\alpha}_1| \leqslant m$ and $|\boldsymbol{\alpha}_2| \leqslant n$.

Let us now define *orthogonality* for anisotropic distributions on $S^*\Omega$ being motivated by the equivalent definition for Radon measures introduced in the previous section.

Anisotropic distributions $\mu$ and $\nu$ in $\mathcal{D}'_{0,\kappa}$ are *orthogonal* ($\mu \perp \nu$) if for any $\varepsilon > 0$ there are open sets $U$ and $V$ in $S^*\Omega$, $U \cup V = S^*\Omega$, such that

$$(\forall \varphi_1, \varphi_2 \in C^{0,\kappa}(S^*\Omega)) \quad \operatorname{supp} \varphi_1 \subseteq U \ \& \ \operatorname{supp} \varphi_2 \subseteq V \implies \begin{aligned} |\langle \mu, \varphi_1 \rangle| &\leqslant \varepsilon \|\varphi_1\|_{C_b^{0,\kappa}(S^*\Omega)} \\ |\langle \nu, \varphi_2 \rangle| &\leqslant \varepsilon \|\varphi_2\|_{C_b^{0,\kappa}(S^*\Omega)} \end{aligned},$$

where $\|\varphi\|_{C_b^{0,\kappa}(S^*\Omega)} = \max_{|\alpha| \leqslant \kappa} \|\partial_\xi^\alpha \varphi(x, \xi)\|_{L^\infty(S^*\Omega)}$ for arbitrary $\varphi \in C^{0,\kappa}(S^*\Omega)$.

Of course, the above definition can easily be extended to include general anisotropic distributions, but we shall use them only in this specific form for H-distributions.

By $\Phi_p \colon L^p_{\operatorname{loc}}(\Omega) \longrightarrow L^{p'}_{\operatorname{loc}}(\Omega)$, for $p \in \langle 1, \infty \rangle$, we denote the Nemyckiĭ operator given by formula $\Phi_p(u) = |u|^{p-2}u$, with appropriate adjustment in the case when $p < 2$ on the set where $u \equiv 0$. Let us emphasise that the Nemyckiĭ operator is a nonlinear continuous operator on these spaces. Furthermore, for arbitrary $u \in L^p_{\operatorname{loc}}(\Omega)$ and $\varphi \in \mathcal{D}(\Omega)$ we have the following bound (see [6, Lemma 4])

$$\|\Phi_p(u)\varphi\|_{L^{p'}(\Omega)} \leqslant \|\Phi_p(u)\varphi\|_{L^p(\Omega)}. \tag{1}$$

**Theorem 2.** *Let $u_n \longrightarrow u$ in $L^p_{\operatorname{loc}}(\Omega)$ and $v_n \longrightarrow v$ in $L^q_{\operatorname{loc}}(\Omega)$ for some $p \in \langle 1, 2 \rangle$ and $q > p'$. Additionally assume $\mu$ is an H-distribution defined by a pair of pure sequences $(u_n - u)$ and $(\Phi_p(u_n - u))$, while $\nu$ is the H-measure defined by pure sequence $(v_n - v)$, and that $\mu \perp \nu$.*

*Then $u_n \bar{v}_n \xrightarrow{*} u \bar{v}$ (vaguely) in the space of Radon measures $C_c(\Omega)'$ (and thus also weakly in the sense of distributions).*

Dem. Consider first the case where $u = v = 0$. Let $\varepsilon > 0$, $\varphi \in C_c(\Omega)$ be nonnegative and let $\{a, b\}$ be a (smooth) partition of unity subordinate to open cover $\{\pi_1(U), \pi_1(V)\}$ of $\Omega$, where $\pi_1 \colon \Omega \times S^{d-1} \longrightarrow \Omega$ is the first projection (to the base of the bundle). Then, as $\mu \perp \nu$, we have (of course, the norm on $C_b(\Omega)$ is just the supremum norm)

$$\begin{aligned} |\langle \mu, \varphi a \boxtimes 1 \rangle| &\leqslant \varepsilon \|\varphi\|_{C_b(\Omega)} \\ |\langle \nu, \varphi b \boxtimes 1 \rangle| &\leqslant \varepsilon \|\varphi\|_{C_b(\Omega)} \end{aligned}.$$

It can easily be seen that

$$\limsup_{n \to \infty} \left| \int_\Omega \varphi u_n \bar{v}_n \, dx \right| = \limsup_{n \to \infty} \left| \int_\Omega \varphi(a + b) u_n \bar{v}_n \, dx \right|$$

$$\leqslant \limsup_{n \to \infty} \left( \int_\Omega \varphi a |u_n| |v_n| \, dx + \int_\Omega \varphi b |u_n| |v_n| \, dx \right) =: \limsup_{n \to \infty} (A_n + B_n)$$

Since $(v_n)$ is bounded in $\mathrm{L}^{p'}_{\mathrm{loc}}(\Omega)$, after applying the Hölder inequality, we have

$$
\begin{aligned}
(2) \quad \limsup_{n\to\infty} A_n &= \limsup_{n\to\infty} \int_\Omega (\varphi a)^{\frac{1}{p}} (\varphi a)^{\frac{1}{p'}} |u_n|\,|v_n|\,dx \\
&\leqslant \limsup_{n\to\infty} \left( \int_\Omega \varphi a |u_n|^p \, dx \right)^{\frac{1}{p}} \left( \int_\Omega \varphi a |v_n|^{p'} \, dx \right)^{\frac{1}{p'}} \\
&\leqslant \lim_{n\to\infty} C \left( \int_\Omega \varphi a |u_n|^p \, dx \right)^{\frac{1}{p}} = C \langle \mu, \varphi a \boxtimes 1 \rangle^{\frac{1}{p}} \leqslant C \varepsilon^{\frac{1}{p}} \|\varphi\|_{\mathrm{C}_b(\Omega)}^{\frac{1}{p}} .
\end{aligned}
$$

Similarly, by using the boundedness of $(u_n)$ in $\mathrm{L}^p_{\mathrm{loc}}(\Omega)$ and the interpolation inequality we get

$$
\begin{aligned}
(3) \quad \limsup_{n\to\infty} B_n &= \limsup_{n\to\infty} \int_\Omega (\varphi b)^{\frac{1}{2}} (\varphi b)^{\frac{1}{2}} |u_n|\,|v_n|\,dx \\
&\leqslant \limsup_{n\to\infty} \left( \int_\Omega (\varphi b)^{\frac{p}{2}} |u_n|^p \, dx \right)^{\frac{1}{p}} \left( \int_\Omega (\varphi b)^{\frac{p'}{2}} |v_n|^{p'} \, dx \right)^{\frac{1}{p'}} \\
&\leqslant \limsup_{n\to\infty} C \left( \int_\Omega \varphi b |v_n|^2 \, dx \right)^{\frac{\theta}{2}} \left( \int_\Omega (\varphi b)^{\frac{q}{2}} |v_n|^q \, dx \right)^{\frac{1-\theta}{q}} \\
&\leqslant \lim_{n\to\infty} C \left( \int_\Omega \varphi b |v_n|^2 \, dx \right)^{\frac{\theta}{2}} = C \langle \nu, \varphi b \boxtimes 1 \rangle^{\frac{\theta}{2}} \leqslant C \varepsilon^{\frac{\theta}{2}} \|\varphi\|_{\mathrm{C}_b(\Omega)}^{\frac{\theta}{2}} ,
\end{aligned}
$$

where $\theta$ is the convexity parameter for $p'$ relative to 2 and $q$, i.e.

$$ \frac{1}{p'} = \frac{\theta}{2} + \frac{1-\theta}{q} \quad . $$

By adding (2) and (3), as $\varepsilon$ can be chosen arbitrary small, we see that for any nonnegative $\varphi \in \mathrm{C}_b(\Omega)$ we have

$$ \lim_{n\to\infty} \int_\Omega \varphi u_n \bar{v}_n \, dx = 0 \quad . $$

However, an arbitrary complex $\varphi \in \mathrm{C}_c(\Omega)$ can be decomposed to its real and imaginary part, $\varphi = \varphi_r + i\varphi_i$, while real functions $\varphi_r$ and $\varphi_i$ can be written as differences of nonnegative functions ($\varphi_r = \varphi_r^+ - \varphi_r^-$ and $\varphi_i = \varphi_i^+ - \varphi_i^-$). Therefore we can conclude that $u_n \bar{v}_n \xrightarrow{*} 0$ (vaguely) in Radon measures.

It remains to prove the claim in the case when either $u$ or $v$ is not identically 0. Take $u \in \mathrm{L}^p_{\mathrm{loc}}(\Omega)$ and $v \in \mathrm{L}^q_{\mathrm{loc}}(\Omega)$ arbitrary, and express the product by using sequences $(u_n - u)$ and $(v_n - v)$. It follows that

$$ \int_\Omega u_n \bar{v}_n \varphi \, dx = \int_\Omega (u_n - u)(\bar{v}_n - \bar{v})\varphi \, dx + \int_\Omega u_n \bar{v} \varphi \, dx + \int_\Omega u \bar{v}_n \varphi - \int_\Omega u \bar{v} \, dx \quad . $$

By the special case, the first term on the right converges to 0. We can easily pass to limits in the second and the third terms, obtaining

$$ u\bar{v} + u\bar{v} - u\bar{v} = u\bar{v} \quad , $$

as claimed.

**Q.E.D.**

Let us remark that for the definition of H-distributions one can take $q = p'$; however, the above proof of the orthogonality result relies on the fact that $q > p'$.

## 4. Localisation principle

Let $I_1 := \mathcal{A}_{|2\pi\xi|^{-1}}$ be the Riesz potential and $R_j := \mathcal{A}_{i\xi_j/|\xi|}$ the $j$-th Riesz transform ($j \in 1..d$). We note that $I_1 \colon \mathrm{L}^q(\mathbf{R}^d) \to \mathrm{L}^{q^*}(\mathbf{R}^d)$ is a continuous operator, for any $q \in \langle 1, d \rangle$ and its conjugate Sobolev exponent $\frac{1}{q^*} = \frac{1}{q} - \frac{1}{d}$. On the other hand, $R_j \colon \mathrm{L}^q(\mathbf{R}^d) \to \mathrm{L}^q(\mathbf{R}^d)$ is continuous for each $q \in \langle 1, \infty \rangle$.

Furthermore, for Schwartz functions $f, g \in \mathcal{S}$ the following identity holds

$$\int_{\mathbf{R}^d} (I_1 f) \, \partial_j g \, d\mathbf{x} = \int_{\mathbf{R}^d} (R_j f) \, g \, d\mathbf{x} \,,$$

and by standard density argument then follows

$$\partial_j I_1(f) = -R_j(f) \,, \qquad \text{for } f \in \mathrm{L}^p(\mathbf{R}^d) \,,$$

so $I_1 \colon \mathrm{L}^p(\mathbf{R}^d) \to \mathrm{W}^{1,p}_{\mathrm{loc}}(\mathbf{R}^d)$ is a continuous operator. The details can be found in [19, chapters 3 and 5].

Many applications of H-measures and H-distributions are based on the so-called localisation principle (cf. [5]), which generalises the compactness by compensation argument to the variable coefficients case. For H-distributions it was first established in a general form in [7, Theorem 4.1], and then written more precisely in [6, Theorem 8], the paper where H-distributions were first recognised as anisotropic distributions in $\mathcal{D}'_{0,Q}(\mathrm{S}^*\Omega)$. We shall need still a slightly different form of the result, which in [6] was not proven, but the reader was instructed to appropriately modify the proof in [7]. For the benefit of the reader we shall provide the complete proof here.

**Lemma 1. (localisation principle)** *Assume that $u_n \longrightarrow 0$ in $\mathrm{L}^p_{\mathrm{loc}}(\Omega)$ and $f_n \longrightarrow 0$ in $\mathrm{W}^{-1,p}_{\mathrm{loc}}(\Omega)$ for some $p \in \langle 1, d \rangle$, and that they satisfy differential relations*

$$\mathsf{a} \cdot \nabla u_n = f_n \,,$$

*where $\mathsf{a} \in \mathrm{C}^1(\Omega; \mathbf{R}^d)$.*

*If $\mu$ is an H-distribution corresponding to some subsequence of $(u_n)$ and $(\Phi_p(u_n))$, then*

$$\operatorname{supp} \mu \subseteq \{ (\mathbf{x}, \boldsymbol{\xi}) \in \mathrm{S}^*\Omega \colon \mathsf{a}(\mathbf{x}) \cdot \boldsymbol{\xi} = 0 \} \,.$$

Dem. Let us take test functions $\varphi_1, \varphi_2 \in \mathrm{C}^\infty_c(\Omega)$ and $\psi \in \mathrm{C}^\kappa(\mathrm{S}^{d-1})$; the latter we then extend to a function homogeneous of order 0 on $\mathbf{R}^d_*$ (which can be achieved by a composition with $\boldsymbol{\xi} \mapsto \boldsymbol{\xi}/|\boldsymbol{\xi}|$). For each $n \in \mathbf{N}$ we define $\Psi_n := \overline{(I_1 \circ \mathcal{A}_{\psi(\frac{\xi}{|\xi|})})(\Phi_p(u_n)\varphi_1)} \in \mathrm{W}^{1,p'}_{\mathrm{loc}}(\Omega)$ and $\Theta_n := \Psi_n \varphi_2 \in \mathrm{W}^{1,p'}_c(\Omega)$.

By taking the duality product of the assumed differential relation by $\Theta_n$ we have

(4)
$$\begin{aligned}
{}_{\mathrm{W}^{-1,p}_{\mathrm{loc}}(\Omega)}\langle f_n, \Theta_n \rangle_{\mathrm{W}^{1,p'}_c(\Omega)} &= {}_{\mathrm{W}^{-1,p}_{\mathrm{loc}}(\Omega;\mathbf{R}^d)}\langle \nabla u_n, \Psi_n \varphi_2 \mathsf{a} \rangle_{\mathrm{W}^{1,p'}_c(\Omega;\mathbf{R}^d)} \\
&= -\int_\Omega u_n \operatorname{div}(\Psi_n \varphi_2 \mathsf{a}) \, d\mathbf{x} \\
&= -\int_\Omega u_n \nabla \Psi_n \cdot (\varphi_2 \mathsf{a}) \, d\mathbf{x} - \int_\Omega u_n \Psi_n \operatorname{div}(\varphi_2 \mathsf{a}) \, d\mathbf{x} \,.
\end{aligned}$$

Our goal is to pass to a limit in (4).

First notice that $(\Psi_n)$ is a bounded sequence in $\mathrm{W}^{1,p'}_{\mathrm{loc}}(\Omega)$ i.e. $(\Theta_n)$ is a bounded sequence in $\mathrm{W}^{1,p'}(\Omega)$. Indeed, since (1) holds, $\Phi_p(u_n)$ is bounded in $\mathrm{L}^{p'}_{\mathrm{loc}}(\Omega)$, so $\Phi_p(u_n)\varphi_1 \in \mathrm{L}^r(\Omega)$ for each $r \in \langle 1, p' ]$. Since $\mathcal{A}_{\psi(\xi/|\xi|)}$ is a continuous operator we have that $\mathcal{A}_{\psi(\xi/|\xi|)}(\Phi_p(u_n)\varphi_1)$ is bounded in $\mathrm{L}^r(\Omega)$ for the same $r$. Let $q^* \geqslant p'$. By continuity of $I_1$, as it was noted above, we see that $(\Psi_n)$ is a bounded sequence in $\mathrm{L}^{q^*}(\Omega)$. Now, it is easily seen that $(\varphi\Psi_n)$ is bounded in $\mathrm{L}^{p'}(\Omega)$ for arbitrary $\varphi \in \mathcal{D}(\Omega)$, since $\varphi$ has compact support. On the other hand, by continuity of the Riesz transform $R_j$, $(\nabla(\varphi\Psi_n))$ is a bounded sequence in $\mathrm{L}^{p'}(\Omega)^d$. Therefore, $(\Psi_n)$ is a bounded sequence in $\mathrm{W}^{1,p'}_{\mathrm{loc}}(\Omega)$ and then has a weakly convergent subsequence in $\mathrm{W}^{1,p'}_{\mathrm{loc}}(\Omega)$. For simplicity, we shall not distinguish in notation the sequence from its subsequence.

Notice that the second term on the right-hand side of (4) goes to 0, since $\operatorname{div}(a\varphi_2)u_n \longrightarrow 0$ in $\mathrm{L}^p(\Omega)$ and we are integrating over a compact subset of $\Omega$.

Furthermore,

$$(5) \qquad \lim_{n\to\infty} {}_{\mathrm{W}^{-1,p}_{\mathrm{loc}}(\Omega)}\langle f_n, \Theta_n\rangle_{\mathrm{W}^{1,p'}_{\mathrm{c}}(\Omega)} = 0.$$

Since,

$${}_{\mathrm{W}^{-1,p}_{\mathrm{loc}}(\Omega)}\langle f_n, \Theta_n\rangle_{\mathrm{W}^{1,p'}_{\mathrm{c}}(\Omega)} = {}_{\mathrm{W}^{-1,p}(\Omega)}\langle \varphi f_n, \Theta_n\rangle_{\mathrm{W}^{1,p'}_0(\Omega)},$$

for arbitrary $\varphi \in \mathcal{D}(\Omega)$ such that $\varphi$ is identically 1 on open set $U$ where $\operatorname{supp}\varphi_2 \subseteq U \subseteq \Omega$ (see [3, Theorem 2]) it can easily be seen that

$$|{}_{\mathrm{W}^{-1,p}_{\mathrm{loc}}(\Omega)}\langle f_n, \Theta_n\rangle_{\mathrm{W}^{1,p'}_{\mathrm{c}}(\Omega)}| \leqslant \|\varphi f_n\|_{\mathrm{W}^{-1,p}(\Omega)}\|\Theta_n\|_{\mathrm{W}^{1,p'}_0(\Omega)}$$
$$\leqslant M\|\varphi f_n\|_{\mathrm{W}^{-1,p}(\Omega)},$$

where $M = \sup_{n\in\mathbf{N}} \|\Theta_n\|_{\mathrm{W}^{1,p'}_0(\Omega)}$. Since $f_n \longrightarrow 0$ in $\mathrm{W}^{-1,p}_{\mathrm{loc}}(\Omega)$ we obtain (6). Without loss of generality we can assume that $(u_n)_n$ and $(\Phi_\mathrm{p}(u_n))_n$ are pure sequences in $\mathrm{L}^p_{\mathrm{loc}}(\Omega)$ and $\mathrm{L}^{p'}_{\mathrm{loc}}(\Omega)$ respectively. Let $\mu$ be a H-distribution defined on sequences $(u_n)$ and $(\Phi_\mathrm{p}(u_n))$. Since,

$$\partial_j \Psi_n = -\overline{\mathcal{A}_{\frac{2\pi i \psi(\frac{\boldsymbol{\xi}}{|\boldsymbol{\xi}|})\xi_j}{|\boldsymbol{\xi}|}}}(\varphi_1 \Phi_\mathrm{p}(u_n)),$$

it is easily seen that

$$\left\langle \mu, \sum_{j=1}^d a_j \bar{\varphi}_1 \varphi_2 \bar{\psi}\xi_j/|\boldsymbol{\xi}|\right\rangle = 0.$$

In other words

$$\operatorname{supp}\mu \subseteq \left\{(\mathbf{x}, \boldsymbol{\xi}) \in \mathrm{S}^*\Omega : \mathsf{a}(\mathbf{x})\cdot\boldsymbol{\xi} = 0\right\}.$$

Q.E.D.

## 5. An application – Homogenisation of heterogenous Boltzmann equations

In this section, we additionally assume that $\Omega$ is connected (so it is a domain, not necessarily bounded, in $\mathbf{R}^d$) and that $1 < p < 2$ is fixed. Furthermore, we take $\mathcal{V}$ to be a nonempty bounded domain in $\mathbf{R}^d$, while by $\mathbf{x} \mapsto \check{\mathbf{x}}$ we denote the standard projection of $\mathbf{R}^d$ onto the $d$-dimensional torus $\mathrm{T}^d = \mathbf{R}^d/\mathbf{Z}^d$. Any function $f$ on $\mathrm{T}^d$ defines the unique 1-periodic function $\tilde{f}$ on $\mathbf{R}^d$ by pull-back $\tilde{f}(\mathbf{x}) := f(\check{\mathbf{x}})$, and vice-versa.

Let us emphasise that the Haar measure on $\mathrm{T}^d$ can be defined as the pushforward of the Lebesgue measure on $[0,1)^d$ by the restriction of the above projection $\check{}$ to $[0,1)^d$, which is then a bijection.

Our aim in this section is to present an application of results of the previous section to a homogenisation problem related to the Boltzmann equation.

Consider the following sequence of first order differential operators on $\Omega \times \mathcal{V}$

$$(6) \qquad L_n = \mathsf{a}(\mathbf{x}, \boldsymbol{\lambda}) \cdot \nabla_{\mathbf{x}} + c\left(\mathbf{x}, \frac{\mathbf{x}}{\varepsilon_n}, \boldsymbol{\lambda}\right),$$

where $\varepsilon_n \longrightarrow 0$ is a sequence in $\mathbf{R}^+$, while $c$ is an operator of the form

$$c(\mathbf{x}, \mathbf{y}, \boldsymbol{\lambda}) = \sigma(\mathbf{x}, \mathbf{y}, \boldsymbol{\lambda}) \, \mathcal{K},$$

where $\sigma \in \mathrm{C}(\Omega \times \mathrm{T}^d \times \mathcal{V})$, $\mathsf{a} \colon \Omega \times \mathcal{V} \longrightarrow \mathbf{R}^d$ is a vector function with analytic components and $\mathcal{K} \colon \mathrm{L}^p(\mathcal{V}) \longrightarrow \mathrm{L}^p(\mathcal{V})$ is a continuous linear operator such that its restriction to $\mathrm{L}^{p'}(\mathcal{V})$ takes values in the same space and

$$(7) \qquad \int_{\mathcal{V}} v \mathcal{K} u \, d\boldsymbol{\lambda} = \int_{\mathcal{V}} u \mathcal{K} v \, d\boldsymbol{\lambda}$$

for arbitrary $u \in \mathrm{L}^p(\mathcal{V})$ and $v \in \mathrm{L}^{p'}(\mathcal{V})$.

Physically, $\mathsf{a}(\mathbf{x}, \lambda)$ represents the local particle velocity (drift), while $\sigma$ is the opacity, representing a measure of how strongly the material resists the passage of radiation (or particles). Therefore, in this section, we assume that $\mathsf{a}$ and $\sigma$ are real functions, although all the results in this section hold in the complex case as well.

In general [24], the homogenisation problem consists in determining an operator $\tilde{L}$ such that if $u_n$ is a solution of

$$L_n u_n = f_n$$

for each $n$, where $f_n \longrightarrow f$ in $\mathrm{W}^{-1,p}_{\mathrm{loc}}(\Omega \times \mathcal{V})$, then the weak limit (or at least a weak accumulation point) $u$ of $u_n$ in $\mathrm{L}^p(\Omega \times \mathcal{V})$ satisfies $\tilde{L}u = f$.

In this case, we shall see that $\tilde{L}$ is also a differential operator of a form similar to $L_n$, i.e.

$$\tilde{L} = \mathsf{a}(\mathbf{x}, \boldsymbol{\lambda}) \cdot \nabla_{\mathbf{x}} + \tilde{c}(\mathbf{x}, \boldsymbol{\lambda}),$$

where $\tilde{c}(\mathbf{x}, \boldsymbol{\lambda}) = \hat{\sigma}(\mathbf{x}, 0, \boldsymbol{\lambda})\mathcal{K}$, $\hat{\phantom{x}}$ denoting the Fourier transform in the second variable $\mathbf{y}$.

The standard interpretation is that a sequence of physical laws given by $L_n u_n = f_n$, after passing to the homogenisation limit, gives the homogenised law $\tilde{L}u = f$. In particular, this encompasses modelling of periodic materials when the period tends to zero.

Under these assumptions, we want to characterise the sequences $L_n$, where $\tilde{L}$ is the homogenised operator, satisfying the following:

$$(8) \qquad \begin{array}{l} \text{For any } (u_n) \text{ such that } u_n \longrightarrow u \text{ in } \mathrm{L}^p(\Omega \times \mathcal{V}) \text{ and } L_n u_n = f_n \longrightarrow f \text{ in } \mathrm{W}^{-1,p}_{\mathrm{loc}}(\Omega \times \mathcal{V}), \\ \text{we have } \tilde{L}u = f. \end{array}$$

Before we state and prove the theorem which gives a necessary condition such that $\tilde{L}u = f$, let us remark on the uniform convergence of Fourier series of infinitely differentiable function $f \in \mathrm{C}^\infty(\mathrm{T}^d)$. Let $\mathrm{S}_N f \in \mathrm{C}^\infty(\mathrm{T}^d)$, where

$$(\mathrm{S}_N f)(\mathbf{x}) := \sum_{\substack{\mathbf{q} \in \mathbf{Z}^d \\ |\mathbf{q}_j| \leq N}} \hat{f}(q) e^{2\pi i \mathbf{q} \cdot \mathbf{x}}$$

for arbitrary $\mathbf{x} \in \mathrm{T}^d$. Of course, $\mathrm{S}_N f \longrightarrow f$ in $\mathrm{L}^2(\mathrm{T}^d)$, so there exists a subsequence $(\mathrm{S}_{N_k} f)$ which converges to $f$ a.e. Since $f \in \mathrm{C}^\infty(\mathrm{T}^d)$, from the bounds on the Fourier coefficients [16, Theorem 3.3.9.], we get that the sequence $(\mathrm{S}_N f)$ is uniformly Cauchy and therefore uniformly convergent to a continuous function, which then implies that $\mathrm{S}_N f \longrightarrow f$ uniformly.

**Theorem 3.** *If a sequence of differential operators $(L_n)$ satisfies* (8), *then*

(9) $\qquad (\forall \mathbf{q} \in \mathbf{Z}^d \setminus \{0\}) \quad \hat{c}(\cdot, \mathbf{q}, \cdot) \not\equiv 0 \implies \mathsf{a} \cdot \mathbf{q} \not\equiv 0\,,$

*where $\hat{c}(\cdot, \mathbf{q}, \cdot) = \hat{\sigma}(\cdot, \mathbf{q}, \cdot)\mathcal{K}$.*

Dem. We shall prove the implication, by proving its contrapositive. In other words, let us assume that there exists $\mathbf{q} \neq 0$ such that $\hat{c}(\cdot, \mathbf{q}, \cdot) \not\equiv 0$ and $\mathsf{a} \cdot \mathbf{q} \equiv 0$. By definition of $c$, then there exists $g \in \mathrm{L}^p(\mathcal{V})$ such that $\hat{\sigma}(\cdot, \mathbf{q}, \cdot)\mathcal{K}g \not\equiv 0$. Starting with this $g$ we can build a sequence of solutions to $L_n u_n = f_n$

$$u_n(\mathbf{x}, \boldsymbol{\lambda}) = e^{-2\pi i \mathbf{q} \cdot \mathbf{x}/\varepsilon_n} g(\boldsymbol{\lambda})\,.$$

Indeed, it easily follows that

$$f_n(\mathbf{x}, \boldsymbol{\lambda}) := L_n u_n(\mathbf{x}, \boldsymbol{\lambda}) = \sigma\left(\mathbf{x}, \frac{\mathbf{x}}{\varepsilon_n}, \boldsymbol{\lambda}\right) e^{-2\pi i \mathbf{q} \cdot \mathbf{x}/\varepsilon_n} \mathcal{K} g\,.$$

By taking the Fourier series expansion with respect to the second variable we have for any $(\mathbf{x}, \mathbf{y}, \boldsymbol{\lambda}) \in \Omega \times \mathrm{T}^d \times \mathcal{V}$ the series of the form

(10) $\qquad \displaystyle\sum_{\mathbf{r} \in \mathbf{Z}^d} \hat{\sigma}(\mathbf{x}, \mathbf{r}, \boldsymbol{\lambda}) e^{2\pi i \mathbf{r} \cdot \mathbf{y}}\,,$

where $\hat{\sigma}(\cdot, \mathbf{r}, \cdot)$ are continuous functions defined on $\Omega \times \mathcal{V}$ for each $\mathbf{r} \in \mathbf{Z}^d$.

Let us first show that $f_n \longrightarrow \hat{\sigma}(\cdot, q, \cdot)\mathcal{K}g$ in $\mathrm{L}^p_{\mathrm{loc}}(\Omega \times \mathcal{V})$ weakly. Due to the density given by Lemma 2 in the Appendix, without loss of generality we can assume that $\sigma = \sigma_1 \boxtimes \sigma_2 \boxtimes \sigma_3 \in \mathcal{D}(\Omega) \boxtimes \mathcal{D}(\mathrm{T}^d) \boxtimes \mathcal{D}(\mathcal{V})$. In our case the series in (10) converges uniformly to $\sigma$ with respect to all three variables, since $\sigma_1$ and $\sigma_3$ are bounded functions. Therefore, without loss of generality we can assume that $\hat{\sigma}(\cdot, \mathbf{r}, \cdot) \equiv 0$ for all but finitely many $\mathbf{r} \in \mathbf{Z}^d$. Since $\mathcal{D}(\Omega) \boxtimes \mathcal{D}(\mathcal{V})$ is dense in $\mathrm{L}^{p'}(\Omega \times \mathcal{V})$ it is enough to take test functions of the form $\varphi(\mathbf{x}, \boldsymbol{\lambda}) = \varphi_1(\mathbf{x})\varphi_2(\boldsymbol{\lambda}) \in \mathcal{D}(\Omega)\boxtimes\mathcal{D}(\mathcal{V})$. Now, it easily follows that

$$\int_{\Omega \times \mathcal{V}} f_n(\mathbf{x}, \boldsymbol{\lambda})\varphi(\mathbf{x}, \boldsymbol{\lambda}) d\mathbf{x} d\boldsymbol{\lambda} = \int_\Omega \sigma_1(\mathbf{x})\sigma_2(\mathbf{x}/\varepsilon_n) e^{-2\pi i \mathbf{q} \cdot \mathbf{x}/\varepsilon_n} \varphi_1(x)\, d\mathbf{x} \int_\mathcal{V} \sigma_3(\boldsymbol{\lambda})\varphi_2(\boldsymbol{\lambda}) \mathcal{K}g(\boldsymbol{\lambda})\, d\boldsymbol{\lambda}$$
$$= \int_\mathcal{V} \sigma_3(\boldsymbol{\lambda})\varphi_2(\boldsymbol{\lambda})\mathcal{K}g(\boldsymbol{\lambda})\, d\boldsymbol{\lambda} \sum_r \int_\Omega \sigma_1(\mathbf{x})\hat{\sigma}_2(\mathbf{r}) e^{2\pi i(\mathbf{r}-\mathbf{q})\cdot\mathbf{x}/\varepsilon_n}\varphi_1(\mathbf{x})\, d\mathbf{x}\,.$$

Notice that the terms of the series are oscillatory integrals of the first kind with phase functions $\phi_\mathbf{r}: \mathbf{R}^d \longrightarrow \mathbf{R}$ depending on $\mathbf{r} \in \mathbf{Z}^d$, i.e. $\phi_\mathbf{r}(\mathbf{x}) = 2\pi(\mathbf{r} - \mathbf{q}) \cdot \mathbf{x}$ for arbitrary $\mathbf{x} \in \mathbf{R}^d$. Therefore, if $\mathbf{r} \neq \mathbf{q}$, by [20, page 341, Prop. 4] we get that

$$\int_\Omega \sigma_1(\mathbf{x})\hat{\sigma}_2(\mathbf{r})\varphi(\mathbf{x}) e^{2\pi i(\mathbf{r}-\mathbf{q})\cdot\mathbf{x}/\varepsilon_n}\, d\mathbf{x} \longrightarrow 0\,.$$

Consequently, we have that $f_n \longrightarrow \hat{\sigma}(\cdot, \mathbf{q}, \cdot)\mathcal{K}g$ in $\mathrm{L}^p_{\mathrm{loc}}(\Omega \times \mathcal{V})$. Clearly, $u_n \longrightarrow 0$ in $\mathrm{L}^p(\Omega \times \mathcal{V})$ by applying the same Proposition from [20], and therefore we get that (8) does not hold.

**Q.E.D.**

It is worth noticing that in the proof of the previous theorem the analyticity of $\mathsf{a}$ was not needed, although for establishing the converse statement it will be essential.

In order to prove the converse statement, we shall use the following stronger condition than (9) on $\sigma$:

(11) $\qquad (\forall \mathbf{q} \in \mathbf{Z}^d \setminus \{0\}) \quad \hat{\sigma}(\cdot, \mathbf{q}, \cdot) \not\equiv 0 \implies \mathsf{a} \cdot \mathbf{q} \not\equiv 0\,.$

This condition excludes the possibility that $\hat{\sigma}(\cdot, \mathbf{q}, \cdot)$ is nonzero but $\hat{c}(\cdot, \mathbf{q}, \cdot)$ vanishes on $\mathrm{L}^p(\mathcal{V})$. However, if we additionally assume that $\mathcal{K}$ is injective on $\mathrm{L}^p(\mathcal{V})$, then the conditions are equivalent. Notice that both (9) and (11) give us conditions under which the principal symbol of $L_n$ or $\tilde{L}$ for fixed $\mathbf{q}$ is not identically zero.

**Theorem 4.** *Suppose (11) holds, then the sequence of differential operators $(L_n)$ satisfies* (8).

Dem. We begin the proof of the converse statement with a simple observation on the H-measure of the sequence $\left(d\left(\mathbf{x}, \frac{\mathbf{x}}{\varepsilon_n}, \boldsymbol{\lambda}\right) - \hat{d}(\mathbf{x}, 0, \boldsymbol{\lambda})\right)$ where $d \in \mathrm{C}(\Omega \times \mathrm{T}^d \times \mathcal{V})$. From [23, Example 2.1], by creating a dummy variable, it is not hard to see that the sequence $\left(d\left(\mathbf{x}, \frac{\mathbf{x}}{\varepsilon_n}, \boldsymbol{\lambda}\right) - \hat{d}(\mathbf{x}, 0, \boldsymbol{\lambda})\right)$ is pure and that its corresponding H-measure is

$$\nu = \sum_{q \neq 0} |\hat{d}(\mathbf{x}, \mathbf{q}, \boldsymbol{\lambda})|^2 \delta_{(\mathbf{q}/|\mathbf{q}|, 0)} \, d\mathbf{x} \, d\boldsymbol{\lambda} \, .$$

Let $(u_n)$ be a sequence such that $u_n \longrightarrow u$ in $\mathrm{L}^p(\Omega \times \mathcal{V})$ and $L_n u_n = f_n$. Notice that the only thing we need to show is that

$$(12) \qquad \int_{\Omega \times \mathcal{V}} \sigma\left(\mathbf{x}, \frac{\mathbf{x}}{\varepsilon_n}, \boldsymbol{\lambda}\right) \mathcal{K} u_n \varphi \, d\mathbf{x} \, d\boldsymbol{\lambda} \longrightarrow \int_{\Omega \times \mathcal{V}} \hat{\sigma}(\mathbf{x}, 0, \boldsymbol{\lambda}) \mathcal{K} u \varphi \, d\mathbf{x} \, d\boldsymbol{\lambda} \, ,$$

for any $\varphi \in \mathcal{D}(\Omega \times \mathcal{V})$.

We divide the proof in two main steps: in the first, we reformulate the problem in a way suitable for applying Theorem 2, while in the second, we show that certain H-distributions are orthogonal in order to apply Theorem 2.

**First step:** Since $\mathcal{V}$ is an open bounded set, for each $N \in \mathbf{N}$ there exists a finite partition $\{K_j\}_{j \in A(N)}$ of $\mathcal{V}$ consisting of Lebesgue measureable sets such that for each $j \in A(N)$ we have that $\operatorname{diam} K_j < \frac{1}{N}$. For convenience in notation, let us define the function $\sigma_j \in \mathrm{C}(\Omega \times \mathrm{T}^d)$, $\sigma_j = \sigma(\cdot, \cdot, \boldsymbol{\lambda}_j)$, where $\boldsymbol{\lambda}_j \in K_j$ for arbitrary $j \in A(N)$. Furthermore, by $\chi_j^N$ we denote the characteristic function of set $K_j$. The main idea of the proof is to choose an appropriate partition $\{K_j\}_{j \in A(N)}$ and prove the statement for functions of the form

$$(13) \qquad (\mathbf{x}, \boldsymbol{\lambda}) \mapsto \sum_{j \in A(N)} \sigma_j\left(\mathbf{x}, \frac{\mathbf{x}}{\varepsilon_n}\right) \chi_j^N(\boldsymbol{\lambda})$$

instead of $\sigma$.

Let $\varphi \in \mathcal{D}(\Omega \times \mathcal{V})$ be arbitrary. Notice that for arbitrary partition $\{K_j\}_{j \in A(N)}$ as above we have that

$$(14) \quad \left| \int_{\Omega \times \mathcal{V}} \sigma\left(\mathbf{x}, \frac{\mathbf{x}}{\varepsilon_n}, \boldsymbol{\lambda}\right) \mathcal{K} u_n \varphi \, d\mathbf{x} d\boldsymbol{\lambda} - \sum_{j \in A(N)} \int_{\Omega \times \mathcal{V}} \sigma_j\left(\mathbf{x}, \frac{\mathbf{x}}{\varepsilon_n}\right) \mathcal{K} u_n \varphi \chi_j^N(\boldsymbol{\lambda}) d\mathbf{x} d\boldsymbol{\lambda} \right|$$
$$\leqslant \int_{\Omega \times \mathcal{V}} \sum_{j \in A(N)} \left| \sigma\left(\mathbf{x}, \frac{\mathbf{x}}{\varepsilon_n}, \boldsymbol{\lambda}\right) - \sigma_j\left(\mathbf{x}, \frac{\mathbf{x}}{\varepsilon_n}\right) \right| \chi_j^N(\boldsymbol{\lambda}) |\mathcal{K} u_n \varphi| \, d\mathbf{x} \, d\boldsymbol{\lambda} \, .$$

Since $\varphi$ has compact support and continuous functions on compact sets are uniformly continuous, by increasing $N$ we can make (14) arbitrarily small. Since $\hat{\sigma}$ is also continuous, without loss of generality we may assume that $\sigma$ is in the form (13).

Furthermore, we may assume that $\varphi = \varphi_1 \boxtimes \varphi_2 \in \mathcal{D}(\Omega) \boxtimes \mathcal{D}(\mathcal{V})$; therefore, since (7) holds we have that

$$\sum_{j \in A(N)} \int_{\Omega \times \mathcal{V}} \sigma_j\left(\mathbf{x}, \frac{\mathbf{x}}{\varepsilon_n}\right) \chi_j^N(\boldsymbol{\lambda}) (\mathcal{K} u_n) \varphi_1 \boxtimes \varphi_2 \, d\mathbf{x} d\boldsymbol{\lambda} = \sum_{j \in A(N)} \int_{\Omega \times \mathcal{V}} \sigma_j\left(\mathbf{x}, \frac{\mathbf{x}}{\varepsilon_n}\right) u_n \varphi_1 \mathcal{K}(\chi_j^N \varphi_2) \, d\mathbf{x} d\boldsymbol{\lambda}.$$

Since $\mathcal{K}(\chi_j^N \varphi_2) \in \mathrm{L}^{p'}(\mathcal{V})$ and $\sigma_j \in \mathrm{C}(\Omega \times \mathrm{T}^d)$ for each $j \in A(N)$, by density, we can assume that $\mathcal{K}(\chi_j^N \varphi_2) \in \mathcal{D}(\mathcal{V})$, which concludes the first step.

**Second step:** Without loss of generality we can take $(u_n - u)$ and $\Phi_p(u_n - u)$ to be pure sequences converging to zero. Let $\mu$ be the H-distribution defined by these two pure sequences. By the localisation principle (Lemma 1) we have that

$$\operatorname{supp} \mu \subseteq \left\{ (\mathbf{x}, \boldsymbol{\lambda}, \boldsymbol{\xi}, \boldsymbol{\eta}) \in \mathrm{S}^*(\Omega \times \mathcal{V}) : \mathsf{a}(\mathbf{x}, \boldsymbol{\lambda}) \cdot \boldsymbol{\xi} = 0 \right\},$$

where $\boldsymbol{\eta}$ is the dual variable to $\boldsymbol{\lambda}$.

For simplicity, let us define functions with a dummy variable, i.e. let $\omega_j \in \mathrm{C}(\Omega \times \mathrm{T}^d \times \mathcal{V})$, $\omega_j(\mathbf{x}, \boldsymbol{\xi}, \boldsymbol{\lambda}) = \sigma_j(\mathbf{x}, \boldsymbol{\xi})$, $j \in A(N)$. Now, take a sequence $v_n^j(\mathbf{x}, \boldsymbol{\lambda}) = \omega_j\left(\mathbf{x}, \frac{\mathbf{x}}{\varepsilon_n}, \boldsymbol{\lambda}\right) - \hat{\omega}_j(\mathbf{x}, 0, \boldsymbol{\lambda})$ in $\mathrm{L}^2_{\mathrm{loc}}(\Omega \times \mathcal{V})$. Clearly from assumption for arbitrary $j \in A(N)$ it follows that,

(15) $\qquad (\forall \mathbf{q} \in \mathbf{Z}^d \setminus \{\mathbf{0}\}) \quad \hat{\omega}_j(\cdot, \mathbf{q}, \cdot) \not\equiv 0 \quad \Longrightarrow \quad \mathsf{a} \cdot \mathbf{q} \not\equiv 0\,.$

By the observation at the beginning, $(v_n^j)$ is a pure sequence for each $j \in A(N)$ and its corresponding H-measure is

$$\nu_j = \sum_{\mathbf{q} \neq 0} |\hat{\omega}_j(\mathbf{x}, \mathbf{q}, \boldsymbol{\lambda})|^2 \delta_{(\mathbf{q}/|\mathbf{q}|, 0)} \, d\mathbf{x} d\boldsymbol{\lambda}\,.$$

It remains to show that sequences $(u_n)$ and $(v_n^j)$ are orthogonal for each $j \in A(N)$. Therefore, let $j \in A(N)$ be arbitrary. For any $\varepsilon > 0$, we define open sets

$$U_\varepsilon = \{(\mathbf{x}, \boldsymbol{\lambda}, \boldsymbol{\xi}, \boldsymbol{\eta}) \in \mathrm{S}^*(\Omega \times \mathcal{V}) : |\mathsf{a}(\mathbf{x}, \boldsymbol{\lambda}) \cdot \boldsymbol{\xi}| > \varepsilon\}$$

and

$$V_\varepsilon = \{(\mathbf{x}, \boldsymbol{\lambda}, \boldsymbol{\xi}, \boldsymbol{\eta}) \in \mathrm{S}^*(\Omega \times \mathcal{V}) : |\mathsf{a}(\mathbf{x}, \boldsymbol{\lambda}) \cdot \boldsymbol{\xi}| < 2\varepsilon\}\,,$$

and also (for given $\mathbf{q} \in \mathbf{Z}^d \setminus \{\mathbf{0}\}$)

$$W_\varepsilon = \left\{(\mathbf{x}, \boldsymbol{\lambda}) \in \Omega \times \mathcal{V} : \left|\mathsf{a}(\mathbf{x}, \boldsymbol{\lambda}) \cdot \frac{\mathbf{q}}{|\mathbf{q}|}\right| < 2\varepsilon \right\}\,.$$

Clearly,

$$|\langle \mu, \varphi \rangle| \leqslant \varepsilon \|\varphi\|_{\mathrm{C}_b(\mathrm{S}^*\Omega)}$$

for $\varphi \in \mathcal{D}(\mathrm{S}^*(\Omega \times \mathcal{V}))$, $\operatorname{supp} \varphi \subseteq U_\varepsilon$.

On the other hand, let $\psi \in \mathcal{D}(\mathrm{S}^*(\Omega \times \mathcal{V}))$, $\operatorname{supp} \psi \subseteq V_\varepsilon$. Now, it easily follows that

$$|\langle \nu, \psi \rangle| \leqslant \sum_{\mathbf{q} \neq 0} \int_{\Omega \times \mathcal{V}} |\hat{\omega}_j(\mathbf{x}, \mathbf{q}, \boldsymbol{\lambda})|^2 |\psi(\mathbf{x}, \boldsymbol{\lambda}, \frac{\mathbf{q}}{|\mathbf{q}|}, 0)| \, d\mathbf{x} d\boldsymbol{\lambda}$$

$$\leqslant \left(\sum_{\mathbf{q} \neq 0} \int_{W_\varepsilon} |\hat{\omega}_j(\mathbf{x}, \mathbf{q}, \boldsymbol{\lambda})|^2 \, d\mathbf{x} d\boldsymbol{\lambda}\right) \|\psi\|_{\mathrm{C}_b(\mathrm{S}^*(\Omega \times \mathcal{V}))}\,.$$

Therefore, it remains to prove that

(16) $\qquad \displaystyle\lim_{\varepsilon \to 0} \sum_{\mathbf{q} \neq 0} \int_{W_\varepsilon} |\hat{\omega}_j(\mathbf{x}, \mathbf{q}, \boldsymbol{\lambda})|^2 \, d\mathbf{x} d\boldsymbol{\lambda} = 0\,.$

Notice that (15) is equivalent to

(17) $\qquad (\forall \mathbf{q} \in \mathbf{Z}^d \setminus \{\mathbf{0}\}) \displaystyle\int_{\{(\mathbf{x}, \boldsymbol{\lambda}) \in \Omega \times \mathcal{V} : \mathsf{a}(\mathbf{x}, \boldsymbol{\lambda}) \cdot \mathbf{q} = 0\}} |\hat{\omega}_j(\mathbf{x}, \mathbf{q}, \boldsymbol{\lambda})|^2 \, d\mathbf{x} d\boldsymbol{\lambda} = 0\,,$

since $\mathsf{a}$ is a vector function with analytic components. As the support of $\varphi_1 \mathcal{K}(\chi_j^N \varphi_2)$ is compact, we may assume that $\omega_j \in \mathrm{L}^2(\Omega \times \mathrm{T}^d \times \mathcal{V})$. Also, for arbitrary $f \in \mathrm{L}^2(\Omega \times \mathrm{T}^d \times \mathcal{V})$ by applying the Plancherel theorem it can easily be seen that

(18) $\qquad \|f\|^2_{\mathrm{L}^2(\Omega \times \mathrm{T}^d \times \mathcal{V})} = \displaystyle\sum_{\mathbf{q} \in \mathbf{Z}^d} \|\hat{f}(\cdot, \mathbf{q}, \cdot)\|^2_{\mathrm{L}^2(\Omega \times \mathcal{V})}.$

Furthermore, because of (17), (18) and the Dominated convergence theorem we have that

$$\lim_{\varepsilon \to 0} \int_{W_\varepsilon} |\hat{\omega}_j(\mathbf{x}, \mathbf{q}, \boldsymbol{\lambda})|^2 \, d\mathbf{x} d\boldsymbol{\lambda} = 0,$$

for arbitrary $\mathbf{q} \in \mathbf{Z}^d \setminus \{\mathbf{0}\}$. Now, it is clear that (16) holds with additional assumption that $\omega_j \in \mathcal{D}(\Omega \times \mathcal{V}) \boxtimes \mathrm{span}(\mathrm{Trig})$. Using a density argument we shall prove the statement for $\omega_j \in \mathrm{L}^2(\Omega \times \mathrm{T}^d \times \mathcal{V})$.

Let $\varepsilon' > 0$ and $h \in \mathcal{D}(\Omega \times \mathcal{V}) \boxtimes \mathrm{span}(\mathrm{Trig})$ (see Appendix 2 for notation) be such that $\|\omega_j - h\|_{\mathrm{L}^2(\Omega \times \mathrm{T}^d \times \mathcal{V})} < \frac{\varepsilon'}{4\|\omega_j\|_{\mathrm{L}^2(\Omega \times \mathrm{T}^d \times \mathcal{V})}}$. Notice that for arbitrary $\mathbf{q} \in \mathbf{Z}^d \setminus \{\mathbf{0}\}$ such that $\hat{h}(\cdot, \mathbf{q}, \cdot) \not\equiv 0$ we may assume that $\hat{\omega}_j(\cdot, \mathbf{q}, \cdot) \not\equiv 0$. Furthermore, let $\delta > 0$ be such that

$$\sum_{\mathbf{q} \neq \mathbf{0}} \int_{W_\varepsilon} |\hat{h}(\mathbf{x}, \mathbf{q}, \boldsymbol{\lambda})|^2 \, d\mathbf{x} d\boldsymbol{\lambda} < \frac{\varepsilon'}{2},$$

for all $0 < \varepsilon < \delta$. Since,

$$\left| \sum_{\mathbf{q} \neq \mathbf{0}} \int_{W_\varepsilon} (|\hat{\omega}_j(\mathbf{x}, \mathbf{q}, \boldsymbol{\lambda})|^2 - |\hat{h}(\mathbf{x}, \mathbf{q}, \boldsymbol{\lambda})|^2) \, d\mathbf{x} d\boldsymbol{\lambda} \right|$$

$$\leqslant \sum_{\mathbf{q} \neq \mathbf{0}} \int_{W_\varepsilon} \left| |\hat{\omega}_j(\mathbf{x}, \mathbf{q}, \boldsymbol{\lambda})| - |\hat{h}(\mathbf{x}, \mathbf{q}, \boldsymbol{\lambda})| \right| \left( |\hat{\omega}_j(\mathbf{x}, \mathbf{q}, \boldsymbol{\lambda})| + |\hat{h}(\mathbf{x}, \mathbf{q}, \boldsymbol{\lambda})| \right) d\mathbf{x} d\boldsymbol{\lambda}$$

$$\leqslant \sum_{\mathbf{q} \neq \mathbf{0}} \int_{W_\varepsilon} \left| \hat{\omega}_j(\mathbf{x}, \mathbf{q}, \boldsymbol{\lambda}) - \hat{h}(\mathbf{x}, \mathbf{q}, \boldsymbol{\lambda}) \right| \left( |\hat{\omega}_j(\mathbf{x}, \mathbf{q}, \boldsymbol{\lambda})| + |\hat{h}(\mathbf{x}, \mathbf{q}, \boldsymbol{\lambda})| \right) d\mathbf{x} d\boldsymbol{\lambda}$$

$$\leqslant \sum_{\mathbf{q} \neq \mathbf{0}} \|\hat{\omega}_j(\cdot, \mathbf{q}, \cdot) - \hat{h}(\cdot, \mathbf{q}, \cdot)\|_{\mathrm{L}^2(\Omega \times \mathcal{V})} \left( \|\hat{\omega}_j(\cdot, \mathbf{q}, \cdot)\|_{\mathrm{L}^2(\Omega \times \mathcal{V})} + \|\hat{h}(\cdot, \mathbf{q}, \cdot)\|_{\mathrm{L}^2(W_\varepsilon)} \right)$$

$$\leqslant \|\omega_j - h\|_{\mathrm{L}^2(\Omega \times \mathrm{T}^d \times \mathcal{V})} \left( \|\omega_j\|_{\mathrm{L}^2(\Omega \times \mathrm{T}^d \times \mathcal{V})} + \sqrt{\sum_{\mathbf{q} \neq \mathbf{0}} \int_{W_\varepsilon} |\hat{h}(\mathbf{x}, \mathbf{q}, \boldsymbol{\lambda})|^2 \, d\mathbf{x} d\boldsymbol{\lambda}} \right) < \frac{\varepsilon'}{2},$$

we have that

$$\sum_{\mathbf{q} \neq \mathbf{0}} \int_{W_\varepsilon} |\hat{\omega}_j(\mathbf{x}, \mathbf{q}, \boldsymbol{\lambda})|^2 \, d\mathbf{x} d\boldsymbol{\lambda} < \varepsilon'.$$

Since $j \in A(N)$ is arbitrary, (12) holds for arbitrary $\varphi \in \mathcal{D}(\Omega \times \mathcal{V})$ by applying Theorem 2, which concludes the proof.

**Q.E.D.**

## Appendix 1: H-distributions

In this Appendix we recall the main results regarding the existence and basic properties of H-distributions, refering the reader to the original paper [6] for further details.

H-distributions can be identified as elements of precise classes of distributions of finite order on the cospherical bundle $S^*\Omega = \Omega \times \mathrm{S}^{d-1}$ (as $\Omega \subseteq \mathbf{R}^d$), so we should first make clear what these classes of *distributions of anisoptropic order* are.

Let $\Omega \subseteq \mathbf{R}^d_\mathbf{x} \times \mathbf{R}^r_\mathbf{y}$ be open. We begin by defining spaces

$$\mathrm{C}^{l,m}(\Omega) := \{ f \in \mathrm{C}(\Omega) : (\forall \boldsymbol{\alpha} \in \mathbf{N}_0^d) \, (\forall \boldsymbol{\beta} \in \mathbf{N}_0^r) \quad \partial_\mathbf{x}^{\boldsymbol{\alpha}} \partial_\mathbf{y}^{\boldsymbol{\beta}} f \in \mathrm{C}(\Omega) \},$$

for arbitrary $l, m \in \mathbf{N}_0$. Each of these spaces can be equipped with its natural Fréchet topology; this is the topology of uniform convergence on compact sets of functions and their derivatives up to order $l$ in $\mathbf{x}$ and $m$ in $\mathbf{y}$.

For a compact set $K \subseteq \Omega$, we define the space $\mathrm{C}_K^{l,m}(\Omega)$, consisting of functions in $\mathrm{C}^{l,m}(\Omega)$ which have support contained in $K$. This subspace inherits the topology from $\mathrm{C}^{l,m}(\Omega)$, which is, when considered only on the subspace, a normable topology, and $\mathrm{C}_K^{l,m}(\Omega)$ is a Banach space.

Of particular importance for defining the notion of H-distribution is the space $\mathrm{C}_c^{l,m}(\Omega)$, consisting of all functions in $\mathrm{C}^{l,m}(\Omega)$ with compact support, which is equipped by a stronger topology than the one inherited from $\mathrm{C}^{l,m}(\Omega)$: the topology of *strict inductive limit*, where we use the topology defined on $\mathrm{C}_K^{l,m}(\Omega)$ for arbitrary nonempty compact set $K \subseteq \Omega$. This space $\mathrm{C}_c^{l,m}(\Omega)$ is also complete.

Any continuous linear functional on $\mathrm{C}_c^{l,m}(\Omega)$ we call a *distribution of anisotropic order*, and such functionals form a vector space $\mathcal{D}'_{l,m}(\Omega) = (\mathrm{C}_c^{l,m}(\Omega))'$.

We can identify any locally summable function $f \in \mathrm{L}^1_{\mathrm{loc}}(\Omega)$ with a distribution of order zero by the usual formula

$$\varphi \mapsto \int_\Omega f(\mathbf{x})\varphi(\mathbf{x})\,d\mathbf{x}\,.$$

It is important to note that this identification relies on the use of the Lebesgue measure on $\Omega$. This leads to continuous and dense embeddings $\mathrm{C}_c^\infty(\Omega) \hookrightarrow \mathrm{C}_c^{l,m}(\Omega) \hookrightarrow \mathcal{D}'(\Omega)$, thus $\mathrm{C}_c^{l,m}(\Omega)$ is a *normal space of distributions* [26], hence its dual $\mathcal{D}'_{l,m}(\Omega)$ forms a subspace of $\mathcal{D}'(\Omega)$, when equipped with the strong topology.

The above construction carries on the manifolds [6, 2.2], and a version of the Schwartz kernel theorem holds [6, Theorem 5]. As a consequence, we have the following theorem on existence of H-distributions [6, Theorem 7] (see also [7] for the first result in this direction).

**Theorem 5.** *If $u_n \longrightarrow 0$ in $\mathrm{L}^p_{\mathrm{loc}}(\mathbf{R}^d)$ and $v_n \stackrel{*}{\rightharpoonup} v$ in $\mathrm{L}^q_{\mathrm{loc}}(\mathbf{R}^d)$ for some $p \in \langle 1,\infty \rangle$ and $q \geqslant p'$, then there exist subsequences $(u_{n'})$, $(v_{n'})$ and a complex valued anisotropic distribution $\mu \in \mathcal{D}'_{0,Q}(\mathrm{S}^*\mathbf{R}^d)$, such that, for any $\varphi_1, \varphi_2 \in \mathrm{C}_c(\mathbf{R}^d)$ and $\psi \in \mathrm{C}^Q(\mathrm{S}^{d-1})$, one has:*

$$\lim_{n'\to\infty} \int_{\mathbf{R}^d} \mathcal{A}_\psi(\varphi_1 u_{n'})(\mathbf{x})\overline{(\varphi_2 v_{n'})(\mathbf{x})}d\mathbf{x} = \lim_{n'\to\infty} \int_{\mathbf{R}^d} (\varphi_1 u_{n'})(\mathbf{x})\overline{\mathcal{A}_{\overline\psi}(\varphi_2 v_{n'})(\mathbf{x})}d\mathbf{x}$$
$$= \langle \mu, \varphi_1 \overline{\varphi}_2 \boxtimes \psi \rangle,$$

*where $\mathcal{A}_\psi : \mathrm{L}^p(\mathbf{R}^d) \longrightarrow \mathrm{L}^p(\mathbf{R}^d)$ is the Fourier multiplier operator with symbol $\psi \in \mathrm{C}^Q(\mathrm{S}^{d-1})$.* ∎

The order $Q = (d-1)(\kappa+2)$ (we assume that $d \geqslant 2$) depends on the smoothness $\kappa := \lfloor d/2 \rfloor + 1$ in the Hőrmander-Mihlin theorem for boundedness of Fourier multipliers [6]. Of course, for $q \in \langle 1,\infty \rangle$, weak and weak-$*$ convergence above coincide since $\mathrm{L}^q_{\mathrm{loc}}(\mathbf{R}^d)$ is reflexive.

If we take $p = q = 2$ in the Theorem, we can also take $v_n = u_n$ and obtain a subsequence such that the above is valid. We can apply the Plancherel theorem and show that $\mu$ is a positive distribution, thus a Radon measure, which is called the *H-measure*. In this case we could have started from a vector sequence $(\mathbf{u}_n)$ (of functions with values in $\mathbf{C}^r$, for details v. [2]), and obtained a hermitian matrix $r \times r$ measure $\boldsymbol{\mu}$.

The Theorem only gives us an upper bound for the order in $\boldsymbol{\xi}$. For example, consider the case $p = q = 2$. In this case, the H-distribution is actually an H-measure, which is of order 0 in $\boldsymbol{\xi}$. As we can see in [6, 5.2], in general we cannot expect to get a distribution of order 0 in $\boldsymbol{\xi}$.

## Appendix 2

In this Appendix we shall state and prove several lemmas which were used in Section 5.

Let $X$ be a smooth manifold. With $\mathrm{C}(X)$ we denote the space of continuous functions $X \longrightarrow \mathbf{C}$. Let $K \subseteq X$ be a nonempty compact subset. By choosing an arbitrary $f \in \mathrm{C}(X)$ we can define

$$p_K(f) := \|f\|_{\mathrm{L}^\infty(K)}.$$

It is not hard to see that $p_K$ is a seminorm on $\mathrm{C}(X)$. Using these seminorms we generate the topology of uniform convergence on compact sets on $\mathrm{C}(X)$.

**Lemma 2.** *Let $X$ and $Y$ be connected smooth manifolds. Then $\mathcal{D}(X) \boxtimes \mathcal{D}(Y)$ is dense in $\mathrm{C}(X \times Y)$ with respect to the topology of uniform convergence on compact sets.*

Dem. It is enough to prove that $\mathcal{D}(X) \boxtimes \mathcal{D}(Y)$ is dense in $\mathcal{D}(X \times Y)$ with respect to the topology of uniform convergence on compact sets. Indeed, let $\varphi \in \mathrm{C}(X \times Y)$, $K \subseteq X \times Y$ compact and $\epsilon > 0$. Since $\mathcal{D}(X \times Y)$ is dense in $\mathrm{C}(X \times Y)$ with respect to the topology of uniform convergence on compact sets (see [6, Theorem 2]) there exists $h \in \mathcal{D}(X \times Y)$ such that $p_K(h - f) < \frac{\epsilon}{2}$. By using the density assumption we can assume there exists $f \in \mathcal{D}(X) \boxtimes \mathcal{D}(Y)$ such that $p_K(f - h) < \frac{\epsilon}{2}$. Now it follows

$$p_K(f - \varphi) \leqslant p_K(f - h) + p_K(h - \varphi) < \epsilon.$$

Therefore, $\mathcal{D}(X) \boxtimes \mathcal{D}(Y)$ is dense in $\mathrm{C}(X \times Y)$.

Now, let us prove that $\mathcal{D}(X) \boxtimes \mathcal{D}(Y)$ is dense in $\mathcal{D}(X \times Y)$ with respect to the topology of uniform convergence on compact sets. Let $(U_i \times V_j, \varphi_i \boxtimes \psi_j)$ be a family of charts of the product manifold $X \times Y$ and let $\{\chi_i \boxtimes \theta_j\}$ be a partition of unity subordinate to $\{U_i \times V_j\}$. Now, let $f \in \mathcal{D}(X \times Y)$ and $K = \operatorname{supp} f$. We see that

$$f|_K = (\sum_{(i,j) \in H} \chi_i \boxtimes \theta_j f)|_K,$$

where $H$ is a finite set. Therefore, we are left with approximating each $\chi_i \boxtimes \theta_j f$. By the Stone-Weierstrass theorem we can approximate each $(\chi_i \boxtimes \theta_j f) \circ (\varphi_i^{-1} \boxtimes \psi_j^{-1})$ with respect to the topology of uniform convergence on $\varphi_i(\operatorname{supp} \chi_i) \times \psi_j(\operatorname{supp} \theta_j)$ with polynomials. Multiplying with adequate bump functions we get the result.

**Q.E.D.**

**Lemma 3.** $\mathcal{D}(\Omega) \boxtimes \mathcal{D}(\mathrm{T}^d)$ *is dense in* $\mathrm{L}^2(\Omega \times \mathrm{T}^d)$.

Dem. It follows directly from the fact that $\mathrm{C}_c(\Omega \times \mathrm{T}^d)$ is dense in $\mathrm{L}^2(\Omega \times \mathrm{T}^d)$ (the measure $\mu$ on $\Omega \times \mathrm{T}^d$ is Radon) and the previous lemma.

Indeed, let $\varepsilon > 0$ and $\varphi \in \mathrm{L}^2(\Omega \times \mathrm{T}^d)$; then there exists $f \in \mathrm{C}_c(\Omega \times \mathrm{T}^d)$ such that $\|\varphi - f\|_{\mathrm{L}^2(\Omega \times \mathrm{T}^d)} < \frac{\varepsilon}{2}$. Furthermore, there exists a sequence $(g_k)$ in $\mathcal{D}(\Omega) \boxtimes \mathcal{D}(\mathrm{T}^d)$ such that $\|g_k - f\|_{\mathrm{L}^\infty(\Omega \times \mathrm{T}^d)} \longrightarrow 0$ and $\operatorname{supp} g_k \subseteq K$ for all $k \in \mathbf{N}$ where $K$ is some compact set contained in $\Omega \times \mathrm{T}^d$. Therefore, there exist $k_0 \in \mathbf{N}$ such that for all $k \in \mathbf{N}, k \geqslant k_0$ we have that $\|f - g_k\|_{\mathrm{L}^\infty(\Omega \times \mathrm{T}^d)} < \frac{\varepsilon}{2\mu(K)}$ i.e. $\|f - g_k\|_{\mathrm{L}^2(\Omega \times \mathrm{T}^d)} < \frac{\varepsilon}{2}$. Therefore,

$$\|\varphi - g_k\|_{\mathrm{L}^2(\Omega \times \mathrm{T}^d)} \leq \|\varphi - f\|_{\mathrm{L}^2(\Omega \times \mathrm{T}^d)} + \|f - g_k\|_{\mathrm{L}^2(\Omega \times \mathrm{T}^d)} < \varepsilon,$$

which proves the result.

**Q.E.D.**

Let $\mathrm{Trig} := \{[0,1\rangle \ni x \mapsto e^{2\pi i q \cdot x} : q \in \mathbf{Z}^d\}$ i.e. Trig is the set of phases of Fourier series. Then, $\mathrm{span}(\mathrm{Trig})$ will, of course, be the linear space of all trigonometric polynomials.

**Lemma 4.** $\mathcal{D}(\Omega) \boxtimes \mathrm{span}(\mathrm{Trig})$ *is dense in the space* $\mathrm{L}^2(\Omega \times \mathrm{T}^d)$.

Dem. Let $\varepsilon > 0$ and $\varphi \in \mathrm{L}^2(\Omega \times \mathrm{T}^d)$. Then, by the previous lemma, there exists $h \in \mathcal{D}(\Omega) \boxtimes \mathrm{L}^2(\mathrm{T}^d)$ such that $\|\varphi - h\|_{\mathrm{L}^2(\Omega \times \mathrm{T}^d)} < \frac{\varepsilon}{2}$. Without loss of generality we can assume $h = f \boxtimes g$, where $f \in \mathcal{D}(\Omega)$ and $g \in \mathrm{L}^2(\mathrm{T}^d)$. Therefore, there exists $N \in \mathbf{N}$ such that $\|g - \mathrm{S}_N g\|_{\mathrm{L}^2(\mathrm{T}^d)} < \frac{\varepsilon}{2\|f\|_{\mathrm{L}^2(\Omega)}}$. Now, it can easily be deduced that $\|h - f \boxtimes \mathrm{S}_N g\|_{\mathrm{L}^2(\Omega \times \mathrm{T}^d)} < \frac{\varepsilon}{2}$. Finally, a simple application of the triangle inequality yields the result.

**Q.E.D.**

**Acknowledgement**

The authors wish to thank Marko Erceg and Aleksandar Bulj for helpful discussions.

**References**

[1] Nenad Antonić: *Memory effects in homogenisation: linear second-order equations*, Arc. Rat. Mech. Analysis **125** (1993) 1–24.

[2] Nenad Antonić: *H-measures applied to symmetric systems*, Proc. Roy. Soc. Edinburgh **126A** (1996) 1133–1155.

[3] Nenad Antonić, Krešimir Burazin: *On certain properties of spaces of locally Sobolev functions*, in *Proceedings of the Conference on applied mathematics and scientific computing*, Zlatko Drmač et al. (eds.), Springer, 2005, pp. 109–120.

[4] Nenad Antonić, Marko Erceg: *One-scale H-distributions and variants*, Results Math. **78**(2023) no. 5, Paper No. 165, 53 pp.

[5] Nenad Antonić, Marko Erceg, Martin Lazar: *Localisation principle for one-scale H-measures*, J. Funct. Anal. **272** (2017) 3410–3454.

[6] Nenad Antonić, Marko Erceg, Marin Mišur: *Distributions of anisotropic order and applications to H-distributions*, Analysis and Applications (2021) 801–843.

[7] Nenad Antonić, Darko Mitrović: *H-distributions: an extension of H-measures to an $L^p - L^q$ setting*, Abs. Appl. Analysis **2011** Article ID 901084 (2011) 12 pp.

[8] Claude Bardos, Harsha Hutridurga: *Simultaneous homogenization and diffusion limits for the linear Boltzmann equation*, Asymptot. Anal. **10** (2016) 111–130.

[9] Noufel Ben Abdallah, Marjolaine Puel, and Michael S. Vogelius: *Diffusion and homogenization limits with separate scales*, Multiscale Model. Simul. **10** (2012) 1148–1179.

[10] Etienne Bernard, Emanuele Caglioti, François Golse: *Homogenization of the linear Boltzmann equation in a domain with a periodic distribution of holes*, SIAM J. Math. Anal. **42** (2010) 2082–2113.

[11] Etienne Bernard, Francesco Salvarani: *Homogenization of the linear Boltzmann equation with a highly oscillating scattering term in extended phase space*, Appl. Math. Lett. **143** (2023) Paper No. 108672, 10 pp.

[12] Nicolas Bourbaki: *Integration I*, Springer, 2004.

[13] Laurent Dumas, François Golse: *Homogenization of transport equations*, SIAM J. Appl. Math., **60** (2000) 1447–1470.

[14] Patrick Gérard: *Microlocal defect measures*, Comm. Partial Diff. Eq. **16** (1991) 1761–1794.

[15] Thierry Goudon, Antoine Mellet: *Homogenization and diffusion asymptotics of the linear Boltzmann equation*, ESAIM Control Optim. Calc. Var. **9** (2003) 371–398.

[16] Loukas Grafakos: *Classical Fourier analysis*, Springer, 2014.

[17] Harsha Hutridurga, Olga Mula, Francesco Salvarani: *Homogenization in the energy variable for a neutron transport model*, Asymptot. Anal. **117** (2019) 1–25.

[18] Laurent Schwartz: *Théorie des distributions*, Hermann, 1966.

[19] Elias M. Stein: *Singular integrals and differentiability properties of functions*, Princeton, 1970.

[20] Elias M. Stein: *Harmonic analysis*, Princeton, 1993.

[21] Luc Tartar: *Nonlocal effects induced by homogenization*, in *Partial differential equations and the calculus of variations* vol. II, pp. 925–938, F. Colombini et al. (eds.), Birkhäuser, 1989.

[22] Luc Tartar: *Memory effects in homogenization*, Arc. Rat. Mech. Analysis **111** (1990) 121–133.

[23] Luc Tartar: *H-measures, a new approach for studying homogenisation, oscillations and concentration effects in partial differential equations*, Proc. Roy. Soc. Edinburgh **115A** (1990) 193–230.

[24] Luc Tartar: *The general theory of homogenization: A personalized introduction*, Springer, 2009.

[25] Luc Tartar: *Approximation of H-measures and beyond*, Ann. Math. Sci. Appl. **8** (2023), no. 3, 461–497.

[26] François Tréves: *Topological vector spaces, distributions and kernels*, Dover, 2006.

[27] Pei Wang, Simon R. Arridge, Ming Jiang: *Radiative transfer equation for media with spatially varying refractive index*, Phys. Rev. **A 90**, 023803 (2014)